# (p, q)-Szász-Mirakyan-Baskakov-Stancu type Operators

Shikha Pandey[a*], Vishnu Narayan Mishra[a,b]

[a] *Department of Applied Mathematics & Humanities, Sardar Vallabhbhai National Institute of Technology, Ichchhanath Mahadev Dumas Road, Surat -395 007 (Gujarat), India.*

[b] *L. 1627 Awadh Puri Colony, Phase –III, Beniganj , Opp. Industrial Training Institute(I.T.I.), Ayodhya Main Road, Faizabad-224 001, (Uttar Pradesh), India.*

**Abstract**

This paper deals with the Stancu variant of (p, q)-Szász-Mirakyan-Baskakov operators. Estimation of moments and establishing few basic approximation results which comprise weighted approximation and direct estimates in view of modulus of continuity is the aim of this paper. Steklov mean method has been used for linear approximation. For the particular case $\alpha = 0$, $\beta = 0$, we get (p, q)-Szász-Mirakyan-Baskakov operators.

## 1. Introduction

In last two decades q-calculus has played very important role in approximation theory. Various operators are introduced in q-form ([2],[4],[7],[9]). (p, q) calculus is the extension of q-calculus. Recently Mursleen et al. depicted (p, q)-analogue of Bernstein Operators in [11]. Elementary notations and operations of (p, q) calculus are:

For $0 < q < p \leq 1$, the (p, q) type integer $[n]_{p,q}$ is described by

$$[n]_{p,q} := \frac{p^n - q^n}{p - q}.$$

(p, q) factorial is expressed as

$$[n]_{p,q}! = [n]_{p,q}[n-1]_{p,q}[n-2]_{p,q} \ldots 1, \quad [0]_{p,q}! = 1.$$

and (p, q) binomial coefficient is expressed as

$$\begin{bmatrix} n \\ k \end{bmatrix}_{p,q} = \frac{[n]_{p,q}!}{[k]_{p,q}![n-k]_{p,q}!}, \quad where \ 0 \leq k \leq n.$$

Expansion of (p, q) power is

$$(x \oplus a)_{p,q}^n = \prod_{s=0}^{n-1}(p^s x + q^s a), \quad and \quad (x \ominus a)_{p,q}^n = \prod_{s=0}^{n-1}(p^s x - q^s a).$$

The definite integral of any function $f$, with limits from 0 to a, is given by

$$\int_0^a f(x) d_{p,q}x = (p-q)a \sum_{k=0}^{\infty} \frac{q^k}{p^{k+1}} f\left(\frac{q^k}{p^{k+1}} a\right), \quad \left|\frac{p}{q}\right| > 1.$$

___________________________

E-mail: sp1486@gmail.com, vishnunarayanmishra@gmail.com

(p, q)-derivative of the function $f$ is given by

$$D_{p,q}f(x) = \frac{f(px) - f(qx)}{(p-q)x}, x \neq 0.$$

Two (p, q)-analogues of the exponential function are, (see [1],[6],[8],[13]):

$$e_{p,q}(x) = \sum_{n=0}^{\infty} \frac{p^{n(n-1)/2}x^n}{[n]_{p,q}}, \quad \text{and} \quad E_{p,q}(x) = \sum_{n=0}^{\infty} \frac{q^{n(n-1)/2}x^n}{[n]_{p,q}}.$$

## 2. Construction of operator

Mohapatra and Walczak studied the Szász-Mirakyan type operators in [10]. Further in [6], (p, q)-Szász-Mirakyan-Baskakov Operators has been illustrated. Motivated by it, we define the stancu variant of (p, q)-Szász-Mirakyan-Baskakov Operators, as

$$T_{n,\alpha,\beta}^{p,q}(f,x) = [n-1]_{p,q} \sum_{n=0}^{\infty} s_{n,k}^{p,q}(x) q^{[k(k+1)-2]/2} p^{(k+1)(k+2)/2} \times \int_0^{\infty} b_{n,k}^{p,q}(t) f\left(\frac{[n]_{p,q}p^k t + \alpha}{[n]_{p,q}+\beta}\right) d_{p,q}t, \tag{1}$$

Where $s_{n,k}^{p,q}(x) = E_{p,q}^{-1}([n]_{p,q}x)\frac{q^{\frac{k(k-1)}{2}}[n]_{p,q}^k x^k}{[k]_{p,q}!}$, $b_{n,k}^{p,q}(t) = \begin{bmatrix} n+k-1 \\ k \end{bmatrix}_{p,q} \frac{t^k}{(1\oplus pt)_{p,q}^{k+n}}$.

For the particular case when $p = q = 1, \alpha = \beta = 0$ of operators (1) is presented in [12].

**Moments**

Using the results of [6], we can easily obtain the first three moments of operators (1), as

**Lemma 1** For $x \in [0, \infty)$ and $0 < q < p \leq 1$ we have

1. $T_{n,\alpha,\beta}^{p,q}(1,x) = 1,$

2. $T_{n,\alpha,\beta}^{p,q}(t,x) = \frac{\alpha}{([n]_{p,q}+\beta)} + \frac{[n]_{p,q}}{qp^2([n]_{p,q}+\beta)[n-2]_{p,q}} + \frac{[n]_{p,q}^2}{pq^2([n]_{p,q}+\beta)[n-2]_{p,q}}x,$

3. $T_{n,\alpha,\beta}^{p,q}(t^2,x) =$

$\frac{\alpha^2}{([n]_{p,q}+\beta)^2} + \frac{2\alpha[n]_{p,q}}{([n]_{p,q}+\beta)^2 qp^2[n-2]_{p,q}} + \frac{[2]_{p,q}[n]_{p,q}^2}{([n]_{p,q}+\beta)^2 p^5 q^3[n-2]_{p,q}[n-3]_{p,q}} + \frac{[n]_{p,q}^2}{([n]_{p,q}+\beta)^2 pq^2[n-2]_{p,q}}\left[\frac{[q(p+[2]_{p,q})+p^2][n]_{p,q}}{p^3q^3[n-3]_{p,q}} + 2\alpha\right]x + \frac{[n]_{p,q}^4}{([n]_{p,q}+\beta)^2 pq^6[n-2]_{p,q}[n-3]_{p,q}}x^2.$

## 3. Direct Result

**Theorem 1** Let $q \in (0,1)$ and $p \in (q, 1]$ the operator maps $C_B$ into $C_B$ and

$$\left\|T_{n,\alpha,\beta}^{p,q}(f)\right\|_{C_B} \leq \|f\|_{C_B},$$

where $C_B[0,\infty)$ is the space of all real valued continuous and bounded functions on $[0,\infty)$ and $\|f\|_{C_B} = \sup_{x \in [0,\infty)}|f(x)|$.

*Proof* From Lemma 1, we have

$$\left|T_{n,\alpha,\beta}^{p,q}(f,x)\right| \leq [n-1]_{p,q} \sum_{k=0}^{\infty} s_{n,k}^{p,q}(x) q^{\frac{[k(k+1)-2]}{2}} p^{\frac{(k+1)(k+2)}{2}} \times \int_0^{\infty} b_{n,k}^{p,q}(t) \left|f\left(\frac{[n]_{p,q}p^k t + \alpha}{[n]_{p,q}+\beta}\right)\right| d_{p,q}t$$

$$\leq \sup_{x \in [0,\infty)}|f(x)| [n-1]_{p,q} \sum_{k=0}^{\infty} s_{n,k}^{p,q}(x) q^{\frac{[k(k+1)-2]}{2}} p^{\frac{(k+1)(k+2)}{2}} \times \int_0^{\infty} b_{n,k}^{p,q}(t) d_{p,q}t$$

$$= \sup_{x \in [0,\infty)}|f(x)| T_{n,\alpha,\beta}^{p,q}(1,x) = \|f\|_{C_B}.$$

□

**Theorem 2** Let $q \in (0,1)$ and $p \in (q,1]$. If $f \in C_B[0,\infty)$, then

$$\left|T_{n,\alpha,\beta}^{p,q}(f,x) - f(x)\right| \le 5\,\omega\left(f, \frac{1}{\sqrt{[n-2]_{p,q}}}\right)\left(\frac{\alpha\sqrt{[n-2]_{p,q}}}{([n]_{p,q}+\beta)} + \frac{[n]_{p,q}}{qp^2([n]_{p,q}+\beta)\sqrt{[n-2]_{p,q}}} + \left|\frac{[n]_{p,q}^2}{pq^2([n]_{p,q}+\beta)\sqrt{[n-2]_{p,q}}} - \sqrt{[n-2]_{p,q}}\right|x\right) +$$

$$\frac{9}{2}\omega_2\left(f, \frac{1}{\sqrt{[n-2]_{p,q}}}\right)\left[\left(\frac{[n]_{p,q}^4}{([n]_{p,q}+\beta)^2 pq^6[n-3]_{p,q}} - \frac{2[n]_{p,q}^2}{pq^2([n]_{p,q}+\beta)} + [n-2]_{p,q}\right)x^2 + \left(\frac{2\alpha[n]_{p,q}^2}{([n]_{p,q}+\beta)^2 pq^2} + \frac{[n]_{p,q}^2[q(p+[2]_{p,q})+p^2]}{([n]_{p,q}+\beta)^2 p^4 q^5[n-3]_{p,q}} - \frac{2\alpha[n-2]_{p,q}}{([n]_{p,q}+\beta)} - \frac{2[n]_{p,q}}{qp^2([n]_{p,q}+\beta)}\right)x + $$

$$\frac{2\alpha[n]_{p,q}}{([n]_{p,q}+\beta)^2 qp^2} + \frac{[2]_{p,q}[n]_{p,q}^2}{([n]_{p,q}+\beta)^2 p^5 q^3[n-3]_{p,q}} + \frac{[n-2]_{p,q}\alpha^2}{([n]_{p,q}+\beta)^2} + 2\right].$$

*Proof* For $x \ge 0$ and $n \in \mathbb{N}$, applying the Steklov mean $f_h$ (See [6]), we have

$$\left|T_{n,\alpha,\beta}^{p,q}(f,x) - f(x)\right| \le T_{n,\alpha,\beta}^{p,q}(|f - f_h|, x) + \left|T_{n,\alpha,\beta}^{p,q}(f_h - f_h(x), x)\right| + |f_h(x) - f(x)|. \tag{3}$$

Using the property of Steklov mean and theorem 1, we get

$$T_{n,\alpha,\beta}^{p,q}(|f - f_h|, x) \le \left\|T_{n,\alpha,\beta}^{p,q}(f - f_h)\right\|_{C_B} \le \|f - f_h\|_{C_B} \le \omega_2(f, h).$$

Now using Taylor's expansion,

$$\left|T_{n,\alpha,\beta}^{p,q}(f_h - f_h(x), x)\right| \le \|f'\|_{C_B} T_{n,\alpha,\beta}^{p,q}(t - x, x) + \frac{1}{2}\|f''\|_{C_B} T_{n,\alpha,\beta}^{p,q}((t - x)^2, x).$$

From the properties of Steklov mean and Lemma 1, we will obtain

$$\left|T_{n,\alpha,\beta}^{p,q}(f_h - f_h(x), x)\right| \le \frac{5}{h}\omega(f, h)\left(\frac{\alpha}{([n]_{p,q}+\beta)} + \frac{[n]_{p,q}}{qp^2([n]_{p,q}+\beta)[n-2]_{p,q}} + \left|\frac{[n]_{p,q}^2}{pq^2([n]_{p,q}+\beta)[n-2]_{p,q}} - 1\right|x\right) + \frac{9}{2h^2}\omega_2(f,h)\,T_{n,\alpha,\beta}^{p,q}((t-x)^2, x),$$

where

$$T_{n,\alpha,\beta}^{p,q}((t - x)^2, x) = \left(\frac{[n]_{p,q}^4}{([n]_{p,q}+\beta)^2 pq^6[n-2]_{p,q}[n-3]_{p,q}} - \frac{2[n]_{p,q}^2}{pq^2([n]_{p,q}+\beta)[n-2]_{p,q}} + 1\right)x^2 + \left(\frac{2\alpha[n]_{p,q}^2}{([n]_{p,q}+\beta)^2 pq^2[n-2]_{p,q}} + \frac{[n]_{p,q}^3[q(p+[2]_{p,q})+p^2]}{([n]_{p,q}+\beta)^2 p^4 q^5[n-2]_{p,q}[n-3]_{p,q}} - \frac{2\alpha}{([n]_{p,q}+\beta)} - \frac{2[n]_{p,q}}{qp^2([n]_{p,q}+\beta)[n-2]_{p,q}}\right)x + \frac{2\alpha[n]_{p,q}}{([n]_{p,q}+\beta)^2 qp^2[n-2]_{p,q}} + \frac{[2]_{p,q}[n]_{p,q}^2}{([n]_{p,q}+\beta)^2 p^5 q^3[n-2]_{p,q}[n-3]_{p,q}} + \frac{\alpha^2}{([n]_{p,q}+\beta)^2}.$$

Replacing $h = \sqrt{\frac{1}{[n-2]_{p,q}}}$ and substituting the above estimated values in (3), we get the above mentioned result. □

### 4. Korovkin Type Weighted Approximation

Define, $\quad C_{x^2}[0,\infty) = \{f : |f(x)| \le M_f(1+x^2),\ M_f > 0, f \text{ is continuous}\},$

and $\quad C_{x^2}^*[0,\infty) = \{f : f \in C_{x^2}[0,\infty), \lim_{|x|\to\infty} \frac{f(x)}{1+x^2} \text{ is finite}\},$

hence for $f \in C_{x^2}^*[0,\infty)$ we have, $\quad \|f\|_{x^2} = \sup_{x \in [0,\infty)} \frac{|f(x)|}{1+x^2}.$

Henceforth, we state the Korovkin type weighted approximation theorem, as in ([3],[5]), for $x \in [0,\infty)$ is,

**Theorem 3** Let $p = (p_n)$ and $q = (q_n)$ satisfies $0 < q_n < p_n \leq 1$ and for $n$ large $p_n \to 1, q_n \to 1$ and $p_n^n \to a, q_n^n \to b$, $\lim_{n\to\infty}[n]_{p_n,q_n} = \infty$. For $f \in C_{x^2}^*[0,\infty)$, we have

$$\lim_{n\to\infty}\left\|T_{n,\alpha,\beta}^{p_n,q_n}(f) - f\right\|_{x^2} = 0.$$

*Proof* It is sufficient to examine the following 3 conditions (See [6])

$$\lim_{n\to\infty}\left\|T_{n,\alpha,\beta}^{p_n,q_n}(t^i) - x^i\right\|_{x^2} = 0, \qquad i = 0,1,2. \tag{2}$$

Since, $T_{n,\alpha,\beta}^{p_n,q_n}(1,x) = 1$, the first condition of (3) is satisfied for $i = 0$. For $n > 3$, we can say

$$\left\|T_{n,\alpha,\beta}^{p_n,q_n}(t) - x\right\|_{x^2} \leq \frac{\alpha}{([n]_{p_n,q_n}+\beta)} + \frac{[n]_{p_n,q_n}}{q_n p_n^2([n]_{p_n,q_n}+\beta)[n-2]_{p_n,q_n}} + \left|\frac{[n]_{p_n,q_n}^2}{p_n q_n^2([n]_{p_n,q_n}+\beta)[n-2]_{p_n,q_n}} - 1\right|x,$$

$$\left\|T_{n,\alpha,\beta}^{p_n,q_n}(t^2) - x^2\right\|_{x^2} \leq$$

$$\frac{1}{([n]_{p_n,q_n}+\beta)^2}\left(\alpha^2 + \frac{2\alpha[n]_{p_n,q_n}}{q_n p_n^2[n-2]_{p_n,q_n}} + \frac{[2]_{p_n,q_n}[n]_{p_n,q_n}^2}{p_n^5 q_n^3[n-2]_{p_n,q_n}[n-3]_{p_n,q_n}}\right) + \frac{[n]_{p_n,q_n}^2}{([n]_{p_n,q_n}+\beta)^2 p_n q_n^2[n-2]_{p_n,q_n}}\left|\frac{[q[n]_{p_n,q_n}(p_n+[2]_{p_n,q_n})+p_n^2]}{p_n^3 q_n^3[n-3]_{p_n,q_n}}\right| +$$

$$2\alpha\left|\sup_{x\in[0,\infty)}\frac{x}{1+x^2} + \left|\frac{[n]_{p_n,q_n}^4}{([n]_{p_n,q_n}+\beta)^2 p_n q_n^6[n-2]_{p_n,q_n}[n-3]_{p_n,q_n}} - 1\right|\frac{x^2}{1+x^2},$$

which results into, for $i = 1,2$ $\qquad \lim_{n\to\infty}\left\|T_{n,\alpha,\beta}^{p_n,q_n}(t^i) - x^i\right\|_{x^2} = 0.$

This completes the proof of Theorem 3. $\qquad\square$

## 5. Voronovskaja-type result

Let $p = (p_n)$ and $q = (q_n)$ satisfies $0 < q_n < p_n \leq 1$ and for $n$ large $p_n \to 1, q_n \to 1$ and $p_n^n \to a, q_n^n \to b$. So, $\lim_{n\to\infty}[n]_{p_n,q_n} \to \infty$ then for any $f \in C_{x^2}^*[0,\infty)$ such that $f', f'' \in C_{x^2}^*[0,\infty)$, we have

$$\lim_{n\to\infty}[n]_{p_n,q_n}\left|T_{n,\alpha,\beta}^{p_n,q_n}(f;x) - f(x)\right| = f'(x)[1 + \alpha + Ax] + \frac{1}{2}f''(x)[1 + Bx]x$$

uniformly on any $[0, K], K > 0$, where

where

$A = \lim_{n\to\infty}[n]_{p_n,q_n}\left(\frac{[n]_{p_n,q_n}}{[n-2]_{p_n,q_n}} - 1\right),$

$B = \lim_{n\to\infty}[n]_{p_n,q_n}\left(\frac{[n]_{p_n,q_n}^3}{([n]_{p_n,q_n}+\beta)[n-2]_{p_n,q_n}[n-3]_{p_n,q_n}} - 1\right).$

*Proof* For $x \in [0,\infty)$, the Taylor's formula for function $f$ is given by

$$f(t) = f(x) + f'(x)(t-x) + \frac{1}{2}f''(x)(t-x)^2 + h(t,x)(t-x)^2,$$

where $h(t,x)$ is Peano form of remainder and $h(\cdot, x) \in C_{x^2}^*[0,\infty)$ and $\lim_{t\to x} h(t,x) = 0$.

Operating the above equation by the operator on both sides we get,

$$[T_{n,\alpha,\beta}^{p_n,q_n}(f;x) - f(x)] = f'(x)T_{n,\alpha,\beta}^{p_n,q_n}(t-x;x) + \frac{1}{2}f''(x)T_{n,\alpha,\beta}^{p_n,q_n}((t-x)^2;x) + T_{n,\alpha,\beta}^{p_n,q_n}(h(t,x)(t-x)^2;x).$$

Using Cauchy-Schwarz inequality we have,

$$T_{n,\alpha,\beta}^{p_n,q_n}(h(t,x)(t-x)^2;x) \leq \sqrt{T_{n,\alpha,\beta}^{p_n,q_n}(h^2(t,x);x)}\sqrt{T_{n,\alpha,\beta}^{p_n,q_n}((t-x)^4;x)}.$$

Since we know $h^2(x,x) = 0$ and $h^2(\cdot,x) \in C_{x^2}^*[0,\infty)$, Hence from theorem (3), $T_{n,\alpha,\beta}^{p_n,q_n}(h^2(t,x);x) = 0$ uniformly for $x \in [0, K]$. So,

$\lim_{n\to\infty} T_{n,\alpha,\beta}^{p_n,q_n}(h(t,x)(t-x)^2;x) = 0.$

Therefore,

$$\lim_{n \to \infty} [n]_{p_n,q_n} \left[ T_{n,\alpha,\beta}^{p_n,q_n}(f;x) - f(x) \right] = f'(x) \lim_{n \to \infty} [n]_{p_n,q_n} T_{n,\alpha,\beta}^{p_n,q_n}(t-x;x) + f''(x) \lim_{n \to \infty} [n]_{p_n,q_n} T_{n,\alpha,\beta}^{p_n,q_n}((t-x)^2;x)$$

$$= f'(x) \lim_{n \to \infty} \left[ \alpha + \frac{[n]_{p_n,q_n}}{q_n p_n^2 [n-2]_{p_n,q_n}} + [n]_{p_n,q_n} \left( \frac{[n]_{p_n,q_n}}{p_n q_n^2 [n-2]_{p_n,q_n}} - 1 \right) x \right]$$

$$+ \frac{1}{2} f''(x) \lim_{n \to \infty} \left[ \frac{1}{([n]_{p_n,q_n} + \beta)} \left( \alpha^2 + \frac{2\alpha [n]_{p_n,q_n}}{q_n p_n^2 [n-2]_{p_n,q_n}} + \frac{(p_n + q_n)[n]_{p_n,q_n}^2}{p_n^5 q_n^3 [n-2]_{p_n,q_n} [n-3]_{p_n,q_n}} \right) \right.$$

$$+ \frac{[n]_{p_n,q_n}^2}{([n]_{p_n,q_n} + \beta) p_n q_n^2 [n-2]_{p_n,q_n}} \left( \frac{[q_n [n]_{p_n,q_n} (2p_n + q_n) + p_n^2]}{p_n^3 q_n^3 [n-3]_{p_n,q_n}} + 2\alpha \right) x - \frac{2[n]_{p_n,q_n}}{q_n p_n^2 [n-2]_{p_n,q_n}} x - 2\alpha x$$

$$+ \left( \frac{[n]_{p_n,q_n}^4}{([n]_{p_n,q_n} + \beta) p_n q_n^6 [n-2]_{p_n,q_n} [n-3]_{p_n,q_n}} + [n]_{p_n,q_n} - \frac{2[n]_{p_n,q_n}^2}{p_n q_n^2 [n-2]_{p_n,q_n}} \right) x^2 \right]$$

$$= f'(x) [1 + \alpha + Ax] + \frac{1}{2} f''(x) [1 + Bx] x,$$

which is the desired result.